\documentclass[11pt]{article}
 \usepackage[ansinew]{inputenc}
 \usepackage[psamsfonts]{amssymb}
 \def\bk{{l\!k}} 

 \def\hH{{H\!H}} 

 \title{Spaces of self-equivalences and   free loops spaces}
 \author{Yves F\'elix and Jean-Claude Thomas}

\begin{document}
 \maketitle

 \begin{abstract} Let  $M$ be a
 simply-connected closed oriented $N$-dimensional  manifold. We prove that
for any field of coefficients there exists a natural
homomorphism of commutative graded  algebras
 $\Psi  : H_\ast (\Omega\, \mbox{aut}_1 M) \to H_{\ast +N}(M^{S^1})$ where
$H_\ast (M^{S^1})$
is the loop algebra  defined by  Chas-Sullivan, \cite{CS}. As usual
$\mbox{aut}_1 X$ (resp. $\Omega X$) denotes  the monoid of the
self-equivalences homotopic to the identity map (resp. the space of based
loops)  of the space $X$. Moreover, if $\bk $ is of
characteristic zero,
$\Psi $ yields isomorphisms $\pi_n(\Omega \mbox{aut}_1 M) \otimes \bk \cong
\hH^{n+N}_{(1)}$ where $ \displaystyle  \oplus
_{l=1}^\infty  \hH^n_{(l)}$ denotes the  Hodge decomposition  on $H^\ast (M
^{S^1})$.
 \end{abstract}

 \vspace{5mm}\noindent {\bf AMS Classification} : 55P35, 55P62, 55P10

 \vspace{2mm}\noindent {\bf Key words} : Free loop space, loop homology,
 self-homotopy equivalences, rational homotopy, Hochschild homology,
$\lambda$-decomposition.

\vspace{5mm}
\noindent{ \bf 1.  Introduction.}  Let $X$ be a path connected space with
base point $x_0$.  We denote   by: $X^{S^1} $ the space of
free loops on $X$, $\Omega X$ the space of based loops of $X$ at $x_0$,
$\mbox{aut}\,  X$ the monoid of self equivalence of $X$
pointed by $Id_X$, $\mbox{aut}_1  X$  the connected component of  $Id_X$ in
$\mbox{aut}  X$,  $\bk$  a field  and $H_\ast(X)$
the singular homology of the space $X$ with coefficients in $\bk$.
Composition of loops or composition of self-equivalences induce
the same  commutative graded algebra structure  on
$H_\ast (\Omega \mbox{aut}_1 X)$.

 Let $M$ be a simply connected $N$-dimensional closed oriented manifold.
It is convenient to  write
$$
\mathbb H_\ast (M) = H_{\ast +N}(M) \mbox{ (resp. } \mathbb H_\ast
(M^{S^1}) = H_{\ast +N}(M^{S^1}) \mbox{)}
$$
Indeed,   $\mathbb H_\ast (M)$  (resp. $\mathbb H_\ast (M^{S^1})$) becomes
a graded commutative  algebra  with  the intersection
product (resp. with the loop product defined by   M. Chas and D. Sullivan,
\cite {CS}).

Our first result reads:

 \vspace{3mm}
\noindent {\bf Theorem 1.} {\sl The natural map
 $$
 \Phi : X \times  \Omega\, \mbox{aut}_1  X  \to
X^{S^1}
\,, \quad (x, \gamma) \mapsto \Phi(x, \gamma) : t\mapsto \gamma(t)(x) \,.
$$
induces   a homomorphism of
commutative graded algebras
$$
 H(\Phi) : \mathbb H_\ast (M)
\otimes
 H_ \ast (\Omega \mbox{\rm aut}_1 M )\to \mathbb H_\ast  (M^{S^1})\,.
$$}

\vspace{3mm}
Let $ 1\in \mathbb H_0(M) =H_d(M)\cong \bk$  be the unit of the algebra
$\mathbb H_\ast(M)$. The homomorphism $H(\Phi)$ restricts
to a homomorphism of commutative  graded algebras
$$
\Psi  : H_\ast (\Omega \mbox{\rm aut}_1 M )\to \mathbb H_\ast (M^{S^1}) \,,
\quad   \Psi(a)=  H(\Phi)(1\otimes a)\,.
$$
Let us denote by,
$$
h : \pi_\ast (\Omega \mbox{aut}_1 M)  \otimes \bk \to  H_\ast (\Omega
\mbox{aut}_1 M )\,,
$$
the homomorphism of (abelian) graded Lie algebras induced  by the Hurewicz
homomorphism.  Thus, composition of the  two
homomorphisms $\Psi $ and $h$  defines  a homomorphism of graded vector spaces:
$$
 \pi_\ast (\Omega \mbox{aut}_1 M) \otimes  \bk \to  H_{\ast+N}(M^{S^1})\,,
$$
which, in turn, induces the dual homomorphism:
$$
 \theta :  H^{\ast +N} (M^{S^1}) \to \left( \pi_\ast (\Omega \mbox{aut}_1 M
\otimes \bk )\right)^\vee \,.
$$

If $\bk$ is a field  of characterictic zero, the Hodge decomposition of the
Hochschild homology
of a commutative graded algebra $A$:
$$
   \hH_\ast (A;A) =\displaystyle  \oplus _{l\geq 0}^\infty  \hH_\ast
^{(l)}(A;A)\,,
$$
was  introduced  by M. Gerstenhaber and S.D. Schack, \cite{G-S}. See also
\cite{Lo} \cite{Lo1}, \cite {Pi} and \cite{Vi1} for
alternative approaches.  From the existence of a commutative model $A$  for
${\cal C}^\ast (M)$ and the isomorphisms:  $\hH
_\ast({\cal C}^\ast (M); {\cal C}^\ast (M)) \cong \hH_\ast (A;A) \cong
H^\ast (M^{S^1})$, it results a Hodge decomposition
$$
\mathbb H^\ast (M ^{S^1}) =  \oplus_{l\geq 0} \mathbb H^\ast _{(l)}(M^{S^1})\,.
$$
 We prove:

 \vspace{3mm}
\noindent {\bf Theorem 2.} { \it If $\bk $ is a field  of characteristic
zero and with notation above, then
$
\theta $
induces  an  isomorphism of graded vector spaces
$$
\left(\pi_\ast (\Omega \mbox{aut}_1 M)\right) \otimes \bk \cong\mathbb
H^{\ast }_{(1)}(M^{S^1})\,, \quad n \geq  0 \,.
$$}

\vspace{3mm}

The remaining of this  paper is organized as follows:

\noindent 2. Proof of theorem 1.

\noindent  3.  Proof of theorem 2.

\noindent  4.  Examples and further comments..

\vspace{5mm}
\noindent {\bf 2.  Proof of theorem 1.}  The proof decomposes in two  steps.

\vspace{3mm}
\underline{First step.}  We denote by $q : M^{S^1} \to M$, the
  free loop space fibration and by $\mbox{Sect}\,(q)$ the space of
    sections of $q$  which are   homotopic to the trivial section $\sigma_0$.
 The  composition  of  loops makes $\mbox{Sect}\,(q)$ into a monoid, with
 multiplication $\mu$ defined by
 $$\mu (\sigma, \tau) (m)(t) =
 \renewcommand{\arraystretch}{1.5}
 \left\{ \begin{array}{ll} \sigma (m)(2t)\,, & t \leq \frac{1}{2}\\
 \tau (m)(2t-1)\,, & t \geq \frac{1}{2}\end{array}\right.
 \renewcommand{\arraystretch}{1}
  \,,\hspace{1cm} \sigma, \tau \in \mbox{Sect}\,(q) \,,  \quad t\in
[0,1]\,, m\in M \,.
$$
Clearly the map
 $$
  \psi : \Omega (\mbox{aut M}_1, id_M) \to \mbox{Sect}\,(q)\,, \quad f
\mapsto \Psi (f) \,, \mbox{ such that } \psi (f)(m)(t) =
f(t)(m)
$$
 is a  homeomorphism of momoids. To prove theorem 1, it suffices to
establishes  that the evaluation map:
$$
ev : M \times \mbox{Sect}\,(q) \to M^{S^1}\,, \quad  (m,\sigma) \mapsto
\sigma (m) \,.
$$
induces a homomorphism of
algebras $ \mathbb H_\ast(M) \otimes  H_\ast(\mbox{Sect}\,(q) )\to \mathbb
H_\ast(M^{S^1})$. This is the purpose of the second
step.

\vspace{3mm}

\underline{Second step.} We recall roughly the definition of the loop
product, by Chas-Sullivan, \cite{CS}, (see also
\cite{CJY} for an another description).

 Let $\alpha  : \triangle ^n \to M^{S^1}$ and
 $\beta : \triangle ^m \to M^{S^1}$ be singular simplices of $M^{S^1}$ and
assume that  the singular simplices of $M$,
  $q\circ \alpha  : \triangle ^n \to M $ and $q\circ \beta  : \triangle ^m
\to M $,  are transverse.  Then  the intersection product  $(q \circ
\alpha) \cdot
(q\circ \beta )$ makes sense and at
 each point
 $(s,t) \in \triangle ^n \times \triangle ^m$ such that  $q\circ \sigma (s)
= q\circ \tau
 (t)$, the composition of the loops  $\alpha  (s) $  can be performed. This
gives a
 chain $\sigma\cdot \tau \in {\cal C}_{m+n-d}(M^{S1} )$ and leads to a
commutative and
  associative multiplication, \cite{CS}-theorem 3.3:
$$
\mathbb H_k(M^{S^1}) \otimes \mathbb H_l (M^{S^1}) \to   \mathbb H_{k+l}
(M^{S^1}) \,, \quad a\otimes b \mapsto a\cdot b\,.
$$

Let $a_1\otimes b_1 \in \mathbb H_k(M)\otimes H_r(\mbox{Sect}(q)) $ and
$a_2\otimes b_2 \in \mathbb H_l(M)\otimes
H_s(\mbox{Sect}(q))
$ then, by definition
$(a_1\otimes b_1)\cdot (a_2\otimes b_2) =(-1) ^{lr} (a_1\cdot a_2 )\otimes
\mu (b_1, b_2)  \in
\mathbb H_{k+l} (M)\otimes H_{r+s}(\mbox{Sect}(q))$ and we have to prove that:
$$
\mathbb H_\ast(ev) \left( (a_1\otimes b_1)\cdot (a_2\otimes b_2)\right) =
\mathbb H_\ast (ev)\left( (a_1\otimes b_1)\right)\cdot
\mathbb H_\ast (ev)\left( (a_1\otimes b_1)\right)
\,.
$$

 First  case: $b_1=b_2=1 \in H_0(\mbox{Sect}(q))$.  M. Chas and D.Sullivan,
\cite{CS}- Proposition 3.4, prove
that
$H_\ast (\sigma_0) :
\mathbb H(M) \to \mathbb H(M^{S^1})$  is
 a homomorphism of graded algebras. Therefore the  restriction of $H_\ast
(ev)$ to $\mathbb H_\ast (M)
 \otimes H_0( \mbox{Sect}\,(q))$  is  a homomorphism of graded algebras.

  Second case: $a_1=a_2=1 \in \mathbb H_0(M)= H_N(M)=\bk \omega $. Let
$\alpha : \triangle ^r \to \mbox{Sect}\,(q)$ and
$\beta  :
\triangle ^s \to
 \mbox{Sect}\,(q)$ be  simplices  and consider the maps  of $M^{S^1}$:
 $$
 \alpha ' : M \times  \Delta^r \stackrel{id \times f}{\longrightarrow} M \times
 \mbox{Sect}\,(q) \stackrel{ev}{\longrightarrow}
  M^{S^1} \,, \quad
  \beta ' : M \to \Delta^s \stackrel{id \times g}{\longrightarrow} M \times
 \mbox{Sect}\,(q)
  \stackrel{ev}{\longrightarrow} M^{S^1}
$$
Since  the simplices,  $q \circ \alpha $ and $q \circ \beta $, are
transverse in $M$ the  product
 $$
f'\cdot g' : M \times \Delta^r\times \Delta^s \stackrel{id \times f \times
 g}{\longrightarrow} M\times \mbox{Sect}\,(q)
 \times \mbox{Sect}\,(q)
 \stackrel{(ev, ev)}{\longrightarrow} M^{S^1}\times_M M^{S^1}
 \stackrel{c}{\longrightarrow} M^{S^1}\,,
$$
is well  defined. Here,  $c$ denotes   pointwise composition of loops.

 As  the multiplication $\mu$ makes commutative the diagram
 $$
\begin{array}{ccc}
 M \times \mbox{Sect}\,(q) \times \mbox{Sect}\,(q) &
 \stackrel{(ev,ev)}{\longrightarrow} & M^{S^1}\times_M M^{S^1}\\
 \downarrow {\scriptstyle id \times \mu} & & \downarrow {\scriptstyle
 c}\\
 M \times \mbox{Sect}\,(q) & \stackrel{ev}{\longrightarrow} & M^{S^1}\,,
 \end{array}
$$
 the map $f'\cdot g'$ is equal to $\mu (f,g)'$.

 Last case:  $a_2=1 \in \mathbb H_0(M)= H_N(M)=\bk \omega $ and $b_1=1 \in
H_0(\mbox{Sect}(q))$.  If  $\alpha :
\triangle ^r \to M$  (resp.  $ \beta  :
\Delta^s \to  \mbox{Sect}\,(q)$) is a simplex of $M$ (resp. of
$\mbox{Sect}\,(q)$)   then the chain  $ev \circ (\alpha \times \beta ) :
 \Delta^r \times \Delta^s \to M^{S^1}$
induced $(ev \otimes (\alpha \times \beta ))'$  coincides with the product of
$ (\sigma _0 \circ \alpha \circ p_2 :  M \times \triangle ^r
\stackrel{p_2}{\to} \triangle ^r \stackrel{ \sigma_0 \circ
\alpha } \to  M^{S^1}$  with
 $\beta  ' : M \times \Delta^s \to M^{S^1}$.  Therefore $H_\ast (ev) ([f]
\otimes [g])
 = H_\ast (\sigma_0)([g]) \cdot H_\ast (ev)(\omega
 \otimes [f])$.  \hfill $\square$

\vspace{5mm}
\noindent{\bf 3. Proof of theorem 2.} Since, $\mathbb Q \subset \bk$, we
may as well suppose that $\bk = \mathbb Q$. Here after we
make    intensive use of the theory  minimal models, in the sense of
Sullivan \cite{Su},  for
 which   refer systematically to \cite{RHT}-$\S 12$. We denote by $(\land
V,d)$ the
 minimal model of $M$. D.
 Sullivan and M.  Vigu\'e, \cite{SV}, have proved that  a relative minimal
model for the fibration
 $q : M^{S^1} \to M$ is given by:
 $$(\land V,d)   \hookrightarrow (\land V \otimes \land sV,D)\,,  \vert
sv\vert = \vert v \vert - 1\,,
 D(v) = d(v)\,, D(sv) = -s(dv)\,,$$
 where $s : \land V \to \land V \otimes \land sV$ is the unique  derivation
defined by: $s(v) = sv$, (see also
\cite{RHT}-$\S 15$(c)). The cochain complex
 $(\land V\otimes \land sV,D)$ decomposes into a direct sum of complexes
 $$(\land V \otimes \land sV,D) = \oplus_{k\geq 0} (\land V \otimes \land^k
 sV,D)\,.$$
 This induces a new  a graduation on $H^\ast (M^{S^1})$, $$H^\ast (M^{S^1}) =
 \oplus_kH^\ast _{(k)}(M^{S^1})\,,
  \hspace{1cm} H^\ast _{(k)}(M^{S^1}) =
  H^\ast (\land V \otimes \land^k sV,D)\,.$$
 We denote by $H_\ast ^{(k)}(M^{S^1})$ the dual graduation on $H_\ast
(M^{S^1})$.
M. Vigu\'e has proved, \cite{Vi1}, that the   Hodge decomposition of the
Hochschild homology
$$
\hH_\ast ( (\land V,d); (\land V, d))= \oplus _{l\geq 0} \hH_\ast^{(l)}
((\land V,d); (\land V, d)\,,
$$
is determined by the isomorphisms   $ H^\ast (\land V \otimes \land^k
sV,D)\cong \hH_\ast ^{(k)} ( (\land V,d); (\land V, d))$.

 By the Milnor-Moore Theorem (\cite{MM}),
 $H_\ast (\mbox{Sect}\,(q);\mathbb Q)$ is isomorphic as an Hopf algebra to the
 universal enveloping algebra on the graded homotopy
  Lie algebra $\pi_\ast (\Omega\, \mbox{aut} M) \otimes \mathbb Q$. Thus
theorem 2 in the introduction is a direct
consequence of theorem 3 below.

 \vspace{3mm}\noindent {\bf Theorem 3.} {\sl The restriction of $H_\ast
(ev)$ to
 $\mathbb H_0(M) \otimes \left(
 \pi_\ast (\mbox{Sect}\,(q)) \otimes \mathbb Q\right)$,
 $$\Phi_1 : \pi_\ast (\mbox{Sect}\,(q)) \otimes \mathbb Q \to \mathbb H_\ast
 (M^{S^1};\mathbb Q)$$
 is an injective homomorphism  whose image is isomorphic to $H_{>
 N}^{(1)}(M^{S^1};\mathbb Q)$.
 }

 \vspace{3mm} \noindent {\bf Proof.} We first construct a quasi-isomorphism
 $\rho : (\land V,d) \to (A,d)$ with $(A,d)$
 a commutative
 differential graded algebra satisfying
 $$\left\{ \begin{array}{l}
 A^0 = \mathbb Q\,, A^1 = 0\,,\\
 A^{>N} = 0\,, A^N = \mathbb Q \Omega\,,\\
 \mbox{dim } A^i < \infty\,, \mbox{ for all } i\,. \end{array} \right.
$$

For this we put,
$\left\{
\begin{array}{ll}
Z^k &= Ker (d : (\land  V)^{k} \to (\land V)^{k+1})\\
 (\land V)^k &= Z^k \oplus S^k\\
I &=  S^{N-1}\oplus dS^{N-1} \oplus S^N \oplus (\land V)^{>N}
\end{array}
\right.
$

Then,  the quotient $(\land V)^N / (S^N \oplus dS^{N-1}) \cong H^N(M) $ has
dimension one. Since  $V^1=0$, the subcomplex   $I$ is
an  ideal of $(\land V,d)$. The acyclicity of $I$ implies that  the natural
projection  $\rho : (\land V,d) \to (A,d)= (\land
V/I,d)$ is a quasi-isomorphism of differential graded algebras.

 The homomorphism  $\rho$ extends to a quasi-isomorphism $\rho \otimes 1 :
(\land V
 \otimes \land sV,D) \to (A \otimes
 \land sV,D)$ with $D(a\otimes sv) = d(a) \otimes sv - (-1)^{\vert a \vert }
 a \cdot (\rho \otimes 1)(sdv)$.
  The complex $(A \otimes \land sV,D)$ also decomposes into the direct sum
 of the complexes $(A \otimes \land^k sV,D)$.

 Denote by $(a_i)\,, i = 1 , \ldots , n$, an homogeneous linear  basis of
$A$, such that
 $a_n = \Omega$, and by $(a_i^\vee)$  the
  dual basis i.e. the linear basis  of $A^\vee = \mbox{Hom}\,(A,\mathbb Q)$
such that $\langle a_i^\vee , a_j\rangle =
\delta_{ij}\,.$

A. Haefliger, \cite{Ha},  has proved that  a model for the evaluation map
$ev : M \times
 \mbox{Sect}\,(q) \to M^{S^1}$ is given by the morphism
 $$
\theta : (A \otimes \land sV,D) \to (A,d) \otimes (\land (A^\vee \otimes
 sV), \delta)\, \quad a\otimes v \mapsto \theta (a\otimes sv) = \sum_i aa_i
\otimes (a_i^\vee\otimes sv)\,.
$$
 Since $D(sV) \subset A \otimes sV$, and $\theta$ is a morphism of
 differential graded algebras,
  $\delta (A^\vee \otimes sV) \subset A^\vee \otimes sV$. We denote by:

 $\rho_1 :
 (\land (A^\vee \otimes sV),\delta ) \to (A^\vee \otimes sV, \delta)$
  the projection onto the complex of indecomposable elements,

 $P : (A,d) \to (\mathbb Q \Omega,0)$
 the homogeneous projection
  onto the component of degree $N$,

 $\pi_1 : (A \otimes\land  sV,D)
 \to (A \otimes sV,D)$ the canonical
 projection on the subcomplex.

 \noindent The dual of $\Phi_1$,
 $$\Phi_1^\vee : H^{\ast +d}(M^{S^1};\mathbb Q ) \to \mbox{Hom}\,
 (\pi_\ast (\mbox{Sect}\,(q)) \otimes \mathbb Q, \mathbb Q )$$
 coincides therefore with $H^\ast (P\circ \rho_1) \circ H^\ast (\theta)$:
 $$(A \otimes \land sV,D) \stackrel{\theta}{\longrightarrow} (A,d) \otimes
 (\land (A^\vee \otimes sV),\delta ) \stackrel{ P \circ
\rho_1}{\longrightarrow}
  \mathbb Q \Omega \otimes (A^\vee \otimes sV,\delta)\,.$$

 \vspace{3mm}\noindent {\bf Lemma.}  {\sl The duality map
$$
Du : A \to A^\vee
 \,, \mbox{\sl such that } \langle Du(a),b\rangle = P(ab) \in \mathbb Q
\Omega \cong \mathbb Q\,,
$$
induces a
 quasi-isomorphism of complexes
 $$Du \otimes 1 : (A \otimes sV,D) \to (A^\vee \otimes sV, \delta)\,.$$
 }

 \vspace{3mm}\noindent {\bf Proof.} Denote by $\alpha_{ij}^k$ and
 $\beta_i^j$   rational numbers defined by the relations
 $\left\{ \begin{array}{l} a_i \cdot a_j = \sum_k \alpha_{ij}^k a_k\\
 d(a_i) = \sum_j \beta_i^j a_j
 \end{array}
 \right.$
 A straighforward computations show that
 \begin{enumerate}
 \item[$\bullet$] $ d(a_i^\vee) = - (-1)^{\vert a_i\vert} \sum_j \beta_j^i\,
 a_j^\vee$
 \item[$\bullet$]  $\sum_r \alpha_{ij}^r \, \alpha_{rk}^t  = \sum_s
 \alpha_{jk}^s\, \alpha_{is}^t\,, $ for   $ i,j,k,t = 1, \ldots , n$
 (associativity)
 \item[$\bullet$]  $\sum_r\alpha_{ij}^r\,\beta_r^s =
 \sum_t\beta_i^t\,\alpha_{tj}^s + (-1)^{\vert a_i\vert}
 \sum_l \beta_j^l\, \alpha_{il}^s\,, $ for  $i,j,l = 1 , \ldots , n$
  (compatibility of the differential  with the multiplication).
 \item[$\bullet$]  $\delta (a_j^\vee \otimes sv) = (-1)^{\vert a_j\vert }
 \left[ \sum_{i,l} \alpha_{il}^j\, (a_l^\vee \otimes sv_i)
  - \sum_r \beta_r^j\, (a_r^\vee \otimes sv)\right]$
 \item[$\bullet$]  $Du(a_i) = \sum_j \alpha_{ij}^n\, a_j^\vee\,.$
 \end{enumerate}

 The duality morphism has degree $N$. A standard computation shows then that
 $$
\delta\circ (Du \otimes 1) =
 (-1)^N (Du \otimes 1) \circ d\,.
$$
 Since  $H^\ast (M)$
 is a Poincar\'e duality algebra and since
 $H^\ast (Du) : H^\ast (M) \to H_\ast (M)$ is the Poincar\'e duality,  $Du
 \otimes 1$ is a quasi-isomorphism. \hfill
  $\square$

 \vspace{3mm}
 It is easy to check  the commutativity of  the
 following diagram of complexes
 $$\renewcommand{\arraystretch}{1.6}
 \begin{array}{ccccc}
 (A\otimes \land sV,D) & \stackrel{\theta}{\longrightarrow} & (A,d)
 \otimes
 (\land (A^{\vee} \otimes sV),
  \delta) & \stackrel{P \otimes \rho_1}{\longrightarrow} & \mathbb Q
 \Omega
 \otimes (A^{\vee} \otimes sV, \delta)\\
 {\scriptstyle \pi_1} \downarrow &&&& \uparrow {\scriptstyle 1
 \otimes (Du
 \otimes 1)}\\
 (A \otimes sV,D) & & \stackrel{\sigma}{\longrightarrow} & &
 \mathbb Q
 \Omega  \otimes (A \otimes sV,D)\,,
 \end{array}
 \renewcommand{\arraystretch}{1}
 $$
 with $\sigma (a\otimes sv) = \Omega \otimes a \otimes sv$. Thus, the
surjectivity of
 $H^\ast ((1 \otimes (Du \otimes 1)) \circ \sigma \circ \pi_1)$  implies
surjectivity of  $\Phi_1^\vee$. Therefore  $\Phi_1$ is an
injective homomorphism whose   image is  isomorphic to
 $H_{>N}^{(1)}(M^{S^1})$. \hfill $\square$

\vspace{5mm}

\noindent{\bf 4.  Examples and further comments.}

 \vspace{5mm} \noindent {\bf Remark 1.}  The morphism $H_\ast (ev)$ is in
general  not
 injective. Let us consider the intersection homomorphism $I : \mathbb
H_\ast (M^{S^1})
 \to H_\ast (\Omega M)$ defines in \cite{CS}-Proposition 3.4. M. Chas and
D. Sullivan have
 proved that $I$ is a morphism of
algebras.
 Together with M. Vigu\'e-Poirrier,\cite{FTV}, we prove that the image of
$I$ is in the center of
 $H_\ast (\Omega M)$ and that the kernel
  of $I$ is a nilpotent ideal  with nilpotency index less than or equal to
$N$. In
 particular, if the center of
 $H_\ast (\Omega M)$ is trivial, then $\mathbb H_{>0} (M^{S^1})$ is a nilpotent
 ideal. On the other hand
  $H_\ast (\mbox{Sect}\,(q))\cong H_\ast(\Omega \mbox{aut}_1M) $ is a free
commutative graded algebra  that we believe generated,
in general,  by
 infinitely many   generators, see  \cite{FT}.

 \vspace{3mm}\noindent {\bf Remark 2.}  In \cite{CJ} Cohen and Jones have
proved
 that $\mathbb H_\ast (M^{S^1})$ is
  isomorphic as an algebra to the Hochschild cohomology $H\!H^\ast ({\cal
 C}^\ast (M), {\cal C}^\ast (M))$.
   On the other hand, in \cite{Ga}, J.B. Gatsinzi establishes, for any space
 $M$ (not necessarily a manifold), an
  algebraic
 isomorphism between $\pi_\ast (\mbox{aut} M, id_M) \otimes \mathbb Q$ and
a subvector
 space of
 $H\!H^ \ast ({\cal C}^\ast (M), {\cal C}^\ast (M))$. The coherence between
the two
 results is given by our Theorem 2.

 \vspace{3mm}\noindent {\bf Problem.} We would like to know if the homomorphism
$$
 H_\ast (\Phi) : \mathbb H_\ast (M)
\otimes
 H_ \ast (\Omega \mbox{\rm aut}_1 M )\to \mathbb H_\ast  (M^{S^1})\,.
$$
 is surjective.

It is the case, for example if $M = \mathbb CP^{N'}\,, N= 2N'$. Other
examples can be checked using for instance \cite{FT}.

 For  spaces, such that $H_\ast (\Phi)$ is surjective, we obtain a
strong connection
 between the behaviour of the sequences of Betti numbers
    dim$\, H_i (M^{S^1})$ and dim$\, \pi_i (\mbox{aut} M, id_M)\otimes
\mathbb Q$.

 \vspace{3mm}\noindent {\bf Example 1.} Let $G$ be a Lie group. The minimal
 model of $G$ is $(\land V,0)$ with $V$
 finite dimensional and concentrated in odd degrees, (Cf. \cite{RHT}-$\S
$12(a)).
 Therefore a model of the free loop space $G^{S^1}$ is $(\land V \otimes
 \land sV, 0)$ and the Haefliger model for
 the space $\mbox{Sect}\, (q)$ is
 $(\land (\land V)^\vee \otimes sV), 0)$. Since the model $\theta$ of the
 evaluation map $ev$ is injective, $H_\ast (ev) : H_\ast (M) \otimes
  H_\ast (\mbox{Sect}\,(q)) \to H_\ast (M^{S^1})$ is surjective. This
implies that
 there exists  an isomorphism of graded algebras.
 $$
\mathbb H_\ast (M^{S^1}) \cong \mathbb H_\ast (M) \otimes H_\ast (\Omega M)\,,
$$
 where the multiplication on the right is simply  the intersection product on
 $\mathbb H_\ast (M)$ and the usual Pontryagin product
 on $H_\ast (\Omega M)$.

 \vspace{3mm}\noindent {\bf Example 2.} Let us assume that  $M$ is a
$\mathbb Q$-hyperbolic space satisfying
either
  $(H^+(M))^3 = 0$ or else  $(H^+(M))^4 = 0$ and $M$ is a coformal space.

 For recall, \cite{RHT}-$\S$35,  a space is $\mathbb Q$-hyperbolic if
dim$\, \pi_\ast (M)\otimes
 \mathbb Q = \infty$. A space  is coformal, \cite{NM},
 if the differential graded algebras ${\cal C}_\ast (\Omega M)$ and
$(H_\ast (\Omega
 M),0)$ are quasi-isomorphic.

Under  the above hypothesis, M. Vigu\'e
 proves, \cite{Vi2}, that
 there exist an integer $n_0$ and
 some constants  $C_1 \geq C_2 > 1$ such that
$$
 C_2^n \leq \sum_{i=1}^n \dim H^i(X^{S^1} ) \leq C_1^n \mbox{  all }  n \geq
 n_0\,.
$$
We directly deduce from Theorem 3 that  the same relations hold for  the
sequence
 dim$\, \pi_i (\mbox{aut} M) \otimes \mathbb Q$ in place of $\dim
H^i(X^{S^1} )$, that is: in both case
the sequences of Betti numbers  have  exponential growth, \cite{RHT}-$\S33$(a).

 \vspace{1cm}

felix@math.ucl.ac.be

 D\'epartement de math\'ematique,

 Universit\'e Catholique de Louvain

 2, Chemin du Cyclotron,

 1348 Louvain-La-Neuve, Belgium

 \vspace{3mm}
jean-claude.thomas@univ-angers.fr

 Facult\'e des Sciences,

 2, Boulevard Lavoisier

 49045 Angers, France


\begin{thebibliography}{99}



 \bibitem{CS} M. Chas and  D. Sullivan,  {\sl String topology}, preprint
  2001


 \bibitem{CJ} R. Cohen and J. Jones,    {\sl A homotopy theoretic
 realization of string topology}, preprint, April  2001.

 \bibitem{CJY} R. Cohen, J. D.S. Jones and J. Yan, {\sl The loop
 homology of spheres and projective spaces}, preprint, February
 2002.



\bibitem{FT} Y. F\'elix and  J.-C. Thomas, {\sl The monoid of
self-equivalences of some homogeneous space}
Expositiones Mathematicae (1994) 305-321.


 \bibitem{RHT} Y. F\'elix, S. Halperin and J.-C. Thomas,  {\it Rational
 Homotopy Theory},
  Graduate Texts in Mathematics 205, Springer-Verlag (2000).

 \bibitem{FTV} Y. F\'elix, J.-C. Thomas and M. Vigu\'e-Poirrier,
  {\sl Structure of the loop homology of a compact manifold}, preprint,
 March 2002.

 \bibitem{Ga} J.-B. Gatsinzi, {\sl The homotopy Lie algebra of classifying
spaces},
 Journal of Pure and Applied Algebra  {\bf 120}  (1997), 281-289.

\bibitem{G-S} M. Gerstenhaber and S.D. Schack, {\it A Hodge type
decomposition  for  commutative  algebras}, Proc. Amer. Math. Soc.
{\bf 98} (1986) 312-316.



  \bibitem{Ha} A. Haefliger, {\sl Rational homotopy of the space of sections of
  a nilpotent bundle},  Trans. Amer. Math. Soc.   {\bf 273} (1982),
 609-620.



\bibitem{Lo} J.-L. Loday,  {\sl Op\'erations sur l'homologie cyclique des
alg\`ebres commutatives},  Invent. Math. {\bf 96}
(1989) 205-230.

\bibitem{Lo1} J.-L. Loday, {\it Cyclic homology}, Grundleren der
mathematischen Wissenschaften {\bf 301} Springer-Verlag 1991.

 \bibitem{MM} J. Milnor and J.C. Moore, {\sl On the structure of Hopf
 algebras},    Annals of Math.  {\bf 81}  (1965)   211-264.

 \bibitem{NM} J. Neisendorfer and T. Miller, {\sl Formal and coformal
spaces}, Illinois J. of Math.
  {\bf 22} (1978), 565-580.

\bibitem{Pi} T. Pirashvili, {\sl Hodge decomposition for higher order
Hochschild homology},   Ann. Sci. \'Ecole Norm. Sup.
 {\bf 33}  (2000)  151-179.

 \bibitem{Su} D. Sullivan, {\sl Infinitesimal computations in topology},
Publ. Math. IHES
 {\bf 47}  (1978), 269-331.

\bibitem{SV} D. Sullivan and M. Vigu\'e, {\sl The homology theory of the
closed geodesic problem},  J. Differential Geom. {\bf
11} (1976) 633-634.

\bibitem{Vi1} M. Vigu\'e-Poirrier, {\sl  D\'ecompositions de l'homologie
cyclique des alg\`
 ebres diff\'erentielles gradu\'ees commutatives},  K-Theory  {\bf 4}
(1991)  399-410.


  \bibitem{Vi2} M. Vigu\'e-Poirrier, {\sl Homotopie rationnelle et
croissance du
 nombre de
 g\'eod\'esiques ferm\'ees},  Ann. Scient. Ec. Norm. Sup.   {\bf 17}
 (1984) 413-431.




  \end{thebibliography}
 \end{document}